\documentclass[12pt]{amsart}
\usepackage{amssymb,amsmath,amsthm}
\newtheorem{lemma}{Lemma}
\def\mn{\medskip\noindent}
\def\sn{\smallskip\noindent}
\def\AF{\mathbb{A}_F}
\def\B{(^*\ ^*_*)}
\def\C{\mathbb{C}}
\def\ClF{\mathrm{Cl}_F}
\def\Coker{\mathrm{Coker}}
\def\Ct{\C^\times}
\def\D{\mathcal{D}}
\def\DEF{\D_{E/F}}
\def\e{\varepsilon}
\def\et{e_1}
\def\F{\mathbb{F}}
\def\Fbar{\overline{\F}}
\def\Ft{\F_2}
\def\Ftbar{\Fbar_2}
\def\Gal{\mathrm{Gal}}
\def\GF{G_F}
\def\GL{\mathrm{GL}}
\def\GLt{\GL_2}
\def\Hom{\mathrm{Hom}}
\def\Im{\mathrm{Im}}
\renewcommand\O{\mathcal{O}}
\def\OF{\O_F}
\def\OFhat{\widehat{\O}_F}
\def\OFt{\O_{F,2}}
\def\OFtpx{\O_{F,\gg 0}^\times}
\def\OKpt{\O_{K',2}}
\def\PGL{\mathrm{PGL}}
\def\PSL{\mathrm{PSL}}
\def\Q{\mathbb{Q}}
\def\Qt{\Q_2}
\def\r{\rho}
\def\tinf{\{2,\infty\}}
\def\U{(^1\ ^*_1)}
\def\vt{v_2}
\def\Z{\mathbb{Z}}
\def\Zt{\Z_2}

\begin{document}
%
%
%
%
\title[Non-existence of mod $2$ Galois representations]{
The non-existence of certain mod $2$ Galois representations
of some small quadratic fields}
\author{Hyunsuk Moon and Yuichiro Taguchi}
\address{Department of Mathematics,
College of Natural Sciences,
Kyungpook National University,
Daegu 702-701, Korea}
\email{hsmoon@knu.ac.kr}
\address{Graduate School of Mathematics,
Kyushu University 33,
Fukuoka 812-8581, Japan}
\email{taguchi@math.kyushu-u.ac.jp}
\subjclass[2000]{11R39, 11F80, 11R32}
\maketitle

\begin{abstract}
For a few quadratic fields, 
the non-existence is proved of continuous irreducible 
mod 2 Galois representations of degree 2 
unramified outside $\{2,\infty\}$.
\end{abstract}

\section{Introduction}

In this paper, 
we prove the following theorem, 
which settles some special cases of versions 
(cf.\ 
Conj.\ 1.1 of \cite{BDJ}, 
Conj.\ 1   of \cite{Schein}  and 
Question 1 in Sect.\ 5 of \cite{Figueiredo}) 
of Serre's modularity conjecture (\cite{S1}, \cite{S2})
for a few quadratic fields:

\mn
{\bf Theorem.}
{\it
Let  $F$  be one of the following quadratic fields:
$$  
\Q(\sqrt{2}),\ 
\Q(\sqrt{3}),\  
\Q(\sqrt{5}),\
\Q(\sqrt{-1}),\ 
\Q(\sqrt{-2}),\ 
\Q(\sqrt{-3}),\ 
\Q(\sqrt{-5}).
$$ 
Then there exist no continuous irreducible representations
$\r:\GF\to\GLt(\Ftbar)$  
unramified outside $\tinf$.
}

\medskip
Here, 
$\GF$  denotes the absolute Galois group  
$\Gal(\overline{F}/F)$  of  $F$, and  
$\Ftbar$  is an algebraic closure of the finite field  
$\F_2$  of two elements. 

The proof is based on the method of discriminant bound 
as in 
\cite{Tate}, \cite{Serre OEuvres}, \cite{B1}, 
\cite{M}, \cite{MT}, \cite{B2}. 
However, 
we need to improve the known upper bounds 
at the prime  $2$. 
This is done in Section 2. 
The proof of the Theorem is given in Section 3. 

It is desirable to have such a theorem 
for mod $p$ representations for other primes  $p$, 
but this seems impossible at least by our method. 

\mn
{\it Acknowledgments.} 
The writing of this paper was greatly
encouraged by 
``Ecole d'\'et\'e sur la conjecture de modularit\'e de Serre",  
Luminy, 9--20, July, 2007. 
We thank its organizers, especially J.-M.\ Fontaine, 
and speakers for their big efforts.
We thank also 
A.\ Ash, 
F.\ Diamond, 
C.\ Khare, 
K.\ Klosin,
T.\ Saito, 
J.-P.\ Serre and 
R.\ Taylor
for their 
information on the literature on 
various versions of Serre's modularity conjecture 
and/or 
comments on the first version of this paper. 

\mn
{\it Convention.} 
For a finite extension  $E/F$  of 
non-Archimedean local fields, 
we denote by  $\DEF$  the different ideal of  $E/F$. 
The $2$-adic valuation  $\vt$
is normalized by  $\vt(2)=1$, and is 
used to measure the order of 
ideals (such as  $\DEF$) in algebraic extensions of 
the $2$-adic field  $\Qt$.   
We denote by
$(^*\ ^*_*)$  and 
$(^1\ ^*_1)$  respectively 
the {\it subgroups} 
$\{(^a_0\ ^b_d)\}$  and  
$\{(^1_0\ ^b_1)\}$  of  $\GLt(\Ftbar)$.

\section{Local lemmas}

Let  $F$  be a finite extension of  $\Qt$,
$D=G_F$  its absolute Galois group, and
$I$      its inertia subgroup.
In this section, we consider mod $2$ representations
$\rho:D\to\GLt(\Ftbar)$
of  $D$.
Let  $E/F$  be the extension cut out by  $\rho$.
We shall estimate the different
$\DEF$  of  $E/F$.
Let  $E_0$  (resp.\ $E_1$)  be the maximal unramified
(resp.\ tamely ramified) subextension of  $E/F$, and
let  $\et=[E_1:E_0]$  be the tame ramification index
of  $E/F$.
Then we have  $\DEF=\D_{E/E_1}\D_{E_1/E_0}$, and
$\vt(\D_{E_1/E_0})=(1-1/\et)/e_F$,
where  $e_F$  is the ramification index of  $F/\Qt$.
Thus it remains for us to calculate  $\D_{E/E_1}$.
We assume  $E/F$  is wildly ramified, with
wild ramification index  $2^m$.
Then the wild inertia subgroup  $G_1$  of  $G:=\Im(\rho)$
is a non-trivial 2-group and,
after conjugation, we may assume
it is contained in  $\U$.
Since  $G_1$  is normal in  $G$  and
the normalizer of  $G_1$  in  $\GLt(\Ftbar)$  is  $\B$,
we may assume that  $\rho$  is of the form
\begin{equation}\label{Eq:rep}
  \rho\ =\ \begin{pmatrix}
              \psi_1 & *      \\
             & \psi_2
         \end{pmatrix},
\end{equation}
where  $\psi_i:D\to\Ftbar^{\times}$
are characters of  $D$. 
Note that the  $\psi_i$'s  have odd order, so that 
they are at most tamely ramified.

\begin{lemma}\label{Lem:wild}
Let the notation be as above. 
Assume further that
$F/\Qt$  has ramification index  $2$. 
If  $E/F$  has ramification index  $2^m$
(i.e.\ if  $e_1=1$), then
there exists a non-negative integer  $m_2\leq m$  such that
\begin{equation*}
  \vt(\DEF)\ =
    \begin{cases}
       \frac{9}{4}-\frac{2^{m_2}+1}{2^m} 
       \quad \text{and}\quad m_2\leq m-1; 
       \text{ or } \\
       2-\frac{2^{m_2}+1}{2^m}.
    \end{cases}
\end{equation*}
If  $\rho$  is non-abelian, 
then the former case does not occur.
\end{lemma}

Here, we say  $\rho$  is (non-)abelian if the group 
$\Im(\rho)$  is (non-)abelian.

\begin{proof} 
By assumption, we have  $E_1=E_0$  and
the characters  $\psi_i$  are unramified.
By local class field theory,
the Galois group  $G_1=\Gal(E/E_1)$, 
which is an elementary $2$-group, is
identified with a quotient of the group
$(1+\pi A)/(1+\pi A)^2$, where  
$A$  is the ring  $\O_{E_1}$  of integers of  $E_1$, 
$\pi$  is a uniformizer of  $A$, and 
$(1+\pi A)^2$  is the subgroup of the 
square elements in the multiplicative group  
$(1+\pi A)$. 
The character group
$X=\Hom(G_1,\Ct)$  of  $G_1$  is identified
with a subgroup of
$\Hom((1+\pi A)/(1+\pi A)^2,\Ct)$.
The subgroup  $X_i$  of  $X$
consisting of the characters with
conductor dividing  $\pi^i$  is identified with
a subgroup of
$\Hom((1+\pi A)^\times/(1+\pi^iA)(1+\pi A)^2,\Ct)$.
It is easy to see that
$$
  \{1\}\ =\ X_1\ \subset\ X_2\ =\ X_3\
                 \subset\ X_4\ \subset\ X_5\ =\ X.
$$
Indeed,
the equality  $X_2=X_3$  follows from the fact that
$(1+\pi^2A)=(1+\pi^3A)(1+\pi A)^2$,
and the equality  $X_5=X$  follows from the fact that
$(1+\pi^5A)\subset(1+\pi A)^2$; 
cf.\ the proof Lemma 2.1 of \cite{M}.
Just as in \cite{Tate}, we can show that
the index  $(X_5:X_4)$  is 1 or 2, 
since the image of  $(1+\pi^4A)$  in
$(1+\pi A)/(1+\pi A)^2$  has order  2.
To see this,
consider the equation
\begin{equation}\label{Eq:sqrt}
  1+a\pi^4\ = \ (1+x\pi^2)^2
\end{equation}
for a given  $a\in A^\times$  and unknown  $x\in A^\times$.
If  $2=c\pi^2$  with  $c\in A^\times$, then
the equation (\ref{Eq:sqrt}) has a solution  $x$
if and only if the congruence
$$
  cx+x^2\ \equiv\ a \pmod{\pi}
$$
has a solution.
Since the $\F_2$-linear map
$\wp:A/\pi A\to A/\pi A$  given by
$x\mapsto cx+x^2$  has
$\dim_{\F_2}\Coker(\wp)=1$,
the equation (\ref{Eq:sqrt}) has a solution for
``half" of the $a$'s.

By assumption, 
$X_5$  has order  $2^m$. Suppose  
$X_2$  has order  $2^{m_2}$.
Then the $2$-adic order of the different  $\DEF=\D_{E/E_1}$
can be calculated as follows by using the
F\"uhrerdiskriminantenproduktformel (\cite{CL}, Chap.\ VI, \S 3):
\\
(1-i)
If  $(X_5:X_4)=2$, then
\begin{align*}
  \vt(\DEF)\ &= \
  \frac{1}{2}\cdot\frac{1}{2^m}
  \left((2^m    -2^{m-1})\times 5 +
        (2^{m-1}-2^{m_2})\times 4 +
        (2^{m_2}-1)      \times 2\right)
  \\
  &= \
  \frac{9}{4}-\frac{2^{m_2}+1}{2^m}.
\end{align*}
(1-ii)
If  $(X_5:X_4)=1$, then
$$
  \vt(\DEF)\ = \
  \frac{1}{2}\cdot\frac{1}{2^m}
  \left((2^m    -2^{m_2})\times 4 +
        (2^{m_2}-1)      \times 2\right)\ = \
  2-\frac{2^{m_2}+1}{2^m}.
$$

Let  $\psi_i$  be the characters in (\ref{Eq:rep}). 
If  $\rho$  is non-abelian
(or equivalently, 
if  $\psi_1\not=\psi_2$  as characters on  $D$), then
$\Gal(E_0/F)=G/G_1$  acts on  $G_1$  (identified with a
subgroup of  $\U$)  via  $\psi_1\psi_2^{-1}$
(cf.\ \cite{OT}, Proof of Prop.\ 2.3).
This induces a similar action on  $X$  which
respects the filtration  $X_i$.
Each orbit by this action has odd cardinality
$|\Im(\psi_1\psi_2^{-1})|$, while
$X_5\smallsetminus X_4$  has $2$-power cardinality 
if it is non-empty.
Thus we must have  $X_5=X_4$, and we are in 
the case (1-ii) above. 
\end{proof}

Specializing the  $F/\Qt$, 
we calculate the value of  $\vt(\DEF)$  more precisely
as follows: 

\begin{lemma}\label{Lem:wild2}
Assume  $F/\Qt$  is a totally ramified quadratic extension.  
Then the extension  $E/F$  has ramification index  $2^m$. 
If  $\rho$  is non-abelian, then there exists a 
non-negative integer  $m_2\leq m$  such that
$$
  \vt(\DEF)\ =\ 2-\frac{2^{m_2}+1}{2^m}.
$$
If  $\rho$  is abelian, then we have  
$m\leq 3$  and  
$\vt(\DEF)\leq 15/8$. In fact, more precisely, we have:
$$
  \vt(\DEF)\ =\
  \begin{cases}
    15/8                      & \text{if  $m=3$},\\
    7/4,\ 3/2\ \text{or}\ 5/4 & \text{if  $m=2$}, \\
    5/4,\ 1\   \text{or}\ 1/2 & \text{if  $m=1$}.
  \end{cases}
$$
\end{lemma}

\begin{proof}
$F/\Qt$  being totally ramified,
any abelian extension of  $F$  has
no non-trivial tame ramification since
$\OF^\times$  is a pro-$2$ group.
Thus the characters  $\psi_i$  in (\ref{Eq:rep}) are
unramified, and  $E/F$  has ramification index  $2^m$.

If  $\rho$  is non-abelian, then 
$\vt(\DEF)$  has the 
second value in Lemma \ref{Lem:wild}.
If  $\rho$  is abelian 
(or equivalently, if  
$\psi_1=\psi_2$  as characters on  $D$), 
then  $G_1$  is
identified with a quotient of
$\OF^\times/(\OF^\times)^2$.
The group  $\OF^\times/(\OF^\times)^2$
has order  $8$.
The different is the largest in the case
where  $G_1\simeq \OF^\times/(\OF^\times)^2$,
in which case  
$m=3$, 
$(X_5:X_4)=(X_4:X_3)=(X_2:X_1)=2$,
and 
$$
  \vt(\DEF)\ =\
  \frac{1}{2}\cdot\frac{1}{8}
  (4\times 5 + 2\times 4 + 1\times 2)\ = \ \frac{15}{8}.
$$
Other cases can be calculated similarly. 
Note that  
$(X_{i+1}:X_i)=1$  or  $2$  
since the residue field of  $\OF$  is  $\Ft$.
\end{proof}

Recall that  
$e_1=[E_1:E_0]$  
denotes the tame ramification index of  $E/F$. 

\begin{lemma}\label{Lem:unr}
If  $F/\Qt$  is the unramified quadratic extension,
then we have  $\et=1$ or $3$.
If  $\rho$  is non-abelian, 
there exist non-negative integers  
$m_2\leq m_4\leq m$  such that
$$
  \vt(\DEF)\ =\
  \begin{cases}
     \ 2 - \frac{1}{2^{m-1}}
    & \text{if  $\et=1$},     \\
          \ \frac{8}{3} - \frac{2^{m_4}+2^{m_2}+1}{3\cdot 2^{m-1}}
    & \text{if  $\et=3$}.
  \end{cases}
$$
If  $\rho$  is abelian, then  $m\leq 3$  and
$\vt(\DEF)\leq 35/12$. 
In fact, more precisely, we have: 
$$
  \vt(\DEF)\ =\
  \begin{cases}
    35/12                   & \text{if  $m=3$}, \\
    8/3   \text{ or } 13/6  & \text{if  $m=2$}, \\
    13/6  \text{ or } 5/3   & \text{if  $m=1$},
  \end{cases}
$$
if  $\et=3$.
If  $\et=1$, then the values of  $\vt(\DEF)$  are
the above values minus  $2/3$.
\end{lemma}

\begin{proof}
By local class field theory,
the characters  $\psi_i$  in  (\ref{Eq:rep})  are
identified with characters of  $F^\times/(1+2\OF)^\times$.
Since  $\OF^\times/(1+2\OF)^\times\simeq\F_4^\times$,
the tamely ramified extension  $E_1/E_0$  has
degree either  1  or  3.

As in the proof of Lemma \ref{Lem:wild},
identify the Galois group  $G_1=\Gal(E/E_1)$  
(resp.\ the character group  $X=\Hom(G_1,\Ct)$)
with a quotient of  $(1+\pi A)/(1+\pi A)^2$ 
(resp.\ a subgroup of  $\Hom((1+\pi A)/(1+\pi A)^2,\Ct)$), 
where  $A=\O_{E_1}$  and
$\pi$  is a uniformizer of  $A$.  
Let  $X_i$  be the subgroup of  $X$  
consisting of the characters of  $G_1$  with
conductor dividing  $\pi^i$.

If  $e_1=1$, then the value of
$\vt(\DEF)$  can be calculated as in Proposition 2.3 of \cite{OT}; 
we have  $\{1\}=X_1\subset X_2\subset X_3=X$  and  
$(X_3:X_2)\leq 2$. 
If  $\rho$  is abelian
(i.e.\ $\psi_1=\psi_2$), then
$X$  is in fact identified with a subgroup of 
the character group of  
$(1+2\OF)/(1+2\OF)^2$, and one has 
$(X_2:X_1)\leq 4$  since  $F$  has residue field  $\F_4$. 
Thus  $|X_3|\leq 8$, and
\begin{equation}\label{Eq:(1)abelian}
  \vt(\DEF)\ =\ 
  \begin{cases}  
    \frac{1}{8}(4\times 3 + 3\times 2)=\frac{9}{4} 
	& \text{if  $m=3$},\\
    \frac{1}{4}(2\times 3+1\times 2)=2 \ \text{ or } & {} \\
    \frac{1}{4}(3\times 2)=\frac{3}{2}
	& \text{if  $m=2$},\\
    \frac{3}{2} \text{ or } \frac{2}{2}=1
	& \text{if  $m=1$}.
  \end{cases}
\end{equation}
If  $\rho$  is non-abelian 
(i.e.\ $\psi_1\not=\psi_2$), 
then as in the last part of the proof
of Lemma \ref{Lem:wild}, we have
$X_3=X_2$, and hence
\begin{equation*}\label{Eq:(1)non-abelian}
  \vt(\DEF)\ 
  =\ \frac{1}{2^m}\left((2^m-1)\times 2\right)\ 
  =\ 2 - \frac{1}{2^{m-1}}.
\end{equation*}

Assume  $e_1=3$.
Then as in the proof of Lemma \ref{Lem:wild}, 
one can show that
$$
  \{1\}\ =\ X_1\ \subset\ X_2\ =\ X_3\
                 \subset\ X_4\ =\ X_5\
                 \subset\ X_6\
                 \subset\ X_7\ =\ X,
$$
with  $(X_7:X_6)=1$  or  $2$. 
By assumption, 
$X_7$  has order  $2^m$. 
Suppose
$|X_2|=2^{m_2}$  and
$|X_4|=2^{m_4}$.
If  $\rho$  is abelian, then
$X$  is identified with a subgroup of 
the character group of  
$(1+2\OF)/(1+2\OF)^2$ 
(so  $X_1=X_2$  and  $X_5=X_6$), and 
$\vt(\D_{E/E_1})$  is calculated to have the
same values as in (\ref{Eq:(1)abelian}).
Adding the tame part
$\vt(\D_{E_1/E_0})=2/3$, we see that  
$\vt(\DEF)$  has the values as in the statement of the lemma.
If  $\rho$  is non-abelian, then as in the former case, we have
$X_7=X_6$, and hence
\begin{align*}
  \vt(\D_{E/E_1})\ &=\
  \frac{1}{3}\cdot\frac{1}{2^m}
  \left((2^m    -2^{m_4})\times 6 +
        (2^{m_4}-2^{m_2})\times 4 +
        (2^{m_2}-1)      \times 2
  \right)
  \\
  &=\ 2 - \frac{2^{m_4}+2^{m_2}+1}{3\cdot 2^{m-1}}.
\end{align*}
Adding the tame part, we obtain
\begin{equation*}
  \vt(\DEF)\
   =\ \frac{8}{3} - \frac{2^{m_4}+2^{m_2}+1}{3\cdot 2^{m-1}}.
\end{equation*}

\end{proof}

\section{Proof of the Theorem}

Suppose there were a continuous irreducible representation
$\rho:\GF\to\GL_2(\Ftbar)$ 
unramified outside  $\tinf$.
Let  $K/F$  be the extension cut out by  $\rho$  and
$G=\Im(\r)$  its Galois group.
As in \cite{Tate},
we distinguish the two cases where
$G$  is solvable and non-solvable.

First we deal with the solvable case.
If  $G$  is solvable, then 
it sits in an exact sequence
$$
   1\to H\to G\to\Z/2\Z\to 1, \qquad
   H\subset \Ftbar^\times\times\Ftbar^\times,
$$
as in Theorem 1 in \S 22 of \cite{Sup}.
Hence  $K$  is an abelian extension of odd degree, 
unramified outside  $\tinf$, 
over the quadratic extension  $K'/F$  corresponding to  $H$.
By using class field theory and noticing that  
$\Q(\sqrt{3})$   has narrow class number 2 (resp.\ 
$\Q(\sqrt{-5})$  has        class number 2), 
we can show that, for each  
$F=\Q(\sqrt{2})$,  $\Q(\sqrt{3})$,  $\Q(\sqrt{5})$  (resp.\ 
$F=\Q(\sqrt{-1})$, $\Q(\sqrt{-2})$, $\Q(\sqrt{-3})$, $\Q(\sqrt{-5})$), 
there are 7 possibilities 
(resp.\   3 possibilities) for such  $K'$.
By examining Jones' tables \cite{JJ}, 
we find them as follows:\\
If  $F=\Q(\sqrt{2})$, then \\
$K'$ =
$\Q(\sqrt{\pm\sqrt{2}})$,  \
$\Q(\sqrt{1\pm\sqrt{2}})$, \
$\Q(\sqrt{2},\sqrt{-1})$,  \
$\Q(\sqrt{-2+\sqrt{2}})$,  \
$\Q(\sqrt{ 2+\sqrt{2}})$; \\
If  $F=\Q(\sqrt{3})$, then \\
$K'$ = 
$\Q(\sqrt{1\pm\sqrt{3}})$,  \
$\Q(\sqrt{-1\pm\sqrt{3}})$, \
$\Q(\sqrt{3},\sqrt{-1})$,   \
$\Q(\sqrt{3},\sqrt{-2})$,   \
$\Q(\sqrt{3},\sqrt{2})$; \\
If  $F=\Q(\sqrt{5})$, then \\
$K'$ = 
$\Q(\sqrt{(1\pm\sqrt{5})/2})$, \ 
$\Q(\sqrt{-1\pm\sqrt{5}})$,    \
$\Q(\sqrt{5},\sqrt{-1})$,      \
$\Q(\sqrt{5},\sqrt{2})$,       \
$\Q(\sqrt{5},\sqrt{-2})$;\\
If  $F=\Q(\sqrt{-1})$, then \\
$K'$ = 
$\Q(\sqrt{-1},\sqrt{2})$, \
$\Q(\sqrt{1\pm\sqrt{-1}})$; \\
If  $F=\Q(\sqrt{-2})$, then \\
$K'$ = 
$\Q(\sqrt{-2},\sqrt{-1})$, \
$\Q(\sqrt{\pm\sqrt{-2}})$; \\
If  $F=\Q(\sqrt{-3})$, then \\
$K'$ = 
$\Q(\sqrt{-3},\sqrt{-1})$, \
$\Q(\sqrt{-3},\sqrt{ 2})$, \
$\Q(\sqrt{-3},\sqrt{-2})$; \\
If  $F=\Q(\sqrt{-5})$, then \\
$K'$ = 
$\Q(\sqrt{-5},\sqrt{-1})$, \
$\Q(\sqrt{-5},\sqrt{ 2})$, \
$\Q(\sqrt{-5},\sqrt{-2})$. \\
All these  $K'$  have class number either 1 or 2. 
Let  $\OKpt=\O_{K'}\otimes_{\Z}\Z_2$
denote the 2-adic completion of the integer ring
$\O_{K'}$  of  $K'$.
Then its multiplicative group
$\OKpt^\times$  is isomorphic to
the direct-product of  $\Zt^{\oplus 4}$  and
a cyclic group of order dividing  12 
(A non-trivial 3-torsion subgroup appears only if   
$K'$  contains  $\Q(\sqrt{-3})$  or  $\Q(\sqrt{5})$). 
Thus there can exist an abelian extension
$K/K'$  of odd degree at most 3. 
But in each case, 
the 3-torsion subgroup of  $\OKpt^\times$  is killed 
(when the reciprocity map is applied) 
by the global unit  
$\zeta_3=(-1+\sqrt{-3})/2$  or  
$\e^2=(3+\sqrt{5})/2$  
(N.B.\ The latter is totally positive). 
Thus there is no abelian extension  $K/K'$  of odd degree
unramified outside  $\tinf$.

Next we prove the non-solvable case.
This is done by the comparison 
of the Tate and Odlyzko bounds for discriminants.
We denote by $d_{K/\Q}$ the discriminant of $K/\Q$, and
$d_K^{1/n} = |d_{K/\Q}|^{1/n}$ the root discriminant of $K$,
where $n=[K: \Q]$. 
By the Odlyzko bound \cite{Odly},
we have
$$
   d_K^{1/n}\ >\
     \begin{cases}
	17.020 & \text{if  $n\geq  120$}, \\
	20.895 & \text{if  $n\geq 1000$}.
     \end{cases}
$$ 
If  $G = \Gal(K/F)$  is non-solvable, 
then  $n=2|G| \geq 120$.
On the other hand, by
Lemmas \ref{Lem:wild2} and \ref{Lem:unr},
we have
$$
 d_K^{1/n}\ \leq\ 
  \begin{cases}
 2\sqrt{2}\cdot 2^2      \ <\ 11.314 & \text{if } F=\Q(\sqrt{2}),  \\
 2\sqrt{3}\cdot 2^2      \ <\ 13.857 & \text{if } F=\Q(\sqrt{3}),  \\
  \sqrt{5}\cdot 2^{35/12}\ <\ 16.885 & \text{if } F=\Q(\sqrt{5}),  \\
 2        \cdot 2^2      \ =\  8     & \text{if } F=\Q(\sqrt{-1}), \\
 2\sqrt{2}\cdot 2^2      \ <\ 11.314 & \text{if } F=\Q(\sqrt{-2}), \\
  \sqrt{3}\cdot 2^{35/12}\ <\ 13.079 & \text{if } F=\Q(\sqrt{-3}), \\
 2\sqrt{5}\cdot 2^2      \ <\ 17.889 & \text{if } F=\Q(\sqrt{-5}).
 \end{cases}
$$
Thus we have a contradiction in all cases but  
$F=\Q(\sqrt{-5})$. 
To deal with the last case, 
let 
$2^m$  be the wild ramification index of  $K/F$  at  2. 
Then the 2-Sylow subgroup of  $G$  has order  $\geq m$. 
If $m\leq 2$, then by Lemma \ref{Lem:wild2} applied to 
a 2-adic completion of  $K/F$, we have 
$\vt(\D_{K/F})\leq 7/4$, and hence 
$$
  d_K^{1/n}\ \leq\ 2\sqrt{5}\cdot 2^{7/4}\ <\ 15.043,
$$
which contradicts the Odlyzko bound.
If  $m\geq 3$, then by 
\S\S 251--253 of \cite{Dickson}, 
the image of  $G$  in  $\PGL_2(\Ftbar)$  contains 
a conjugate of  $\PSL_2(\F_8)$, which has order  504. 
Hence the Odlyzko bound applies with  $n=2|G|>1000$, whence 
a contradiction in the remaining case as well.
\qed


\end{document}